\documentstyle[amsfonts,12pt]{article}
\newcommand{\f}{{\cal O}}
\newcommand{\lon}{\mbox{$\longrightarrow$}}
\newcommand{\rar}{\rightarrow}
\newcommand{\hook}{\hookrightarrow}

\newcommand{\OM}{\Omega^1 M}
\newcommand{\rank}{\mbox{rank}}

\newcommand{\bnu}{\nu}
\newcommand{\bmu}{\mu}

\newtheorem{theorem}{Theorem}

\begin{document}

\title{Twistor theory, complex homogeneous \\manifolds and $G$-structures}

\author{Sergey A.\ Merkulov\\
\small   School of Mathematics and Statistics, University of Plymouth \\
\small                Plymouth, Devon PL4 8AA, United Kingdom
}

\date{ }
\maketitle
\sloppy

\paragraph{0.\ Introduction.}
One of the most useful characteristics of an affine connection on a manifold
$M$ is its (restricted) holonomy group which is defined, up to a conjugation,
as a subgroup of $GL(T_t M)$ consisting of all automorphisms of the tangent
space $T_t M$ at a point $t\in M$ induced by parallel translations along the
$t$-based contractible loops in $M$. Which groups can occur as holonomies of
affine connections? By Hano and Ozeki \cite{HO}, any closed subgroup of a
general linear group can be realised as a holonomy of some affine connection
(which in general has a non-vanishing torsion tensor). The same question, if
posed in the class of {\em torsion-free}\, (non-locally symmetric) affine
connections only, is
not yet answered. According to Berger \cite{Berger}, the list of all possible
irreducibly acting holonomies of such connections is very
restricted. How much is known about this list?
In his seminal paper \cite{Berger}, Berger found a list of groups which
embraces all possible holonomies of torsion-free {\em metric}\, connections,
though his approach provides no method to distinguish which entries
can indeed be realised as holonomies and which are superfluous.
Later much work has been done to refine this list and to prove existence
of Riemannian metrics with special holonomies \cite{Al,Bryant1,Bryant2,S}.
In the same paper Berger presented
also a list of all but a finite number of possible candidates to irreducible
holonomies of "non-metric" torsion-free affine connections. How many holonomies
are missing from this second list is not known, but, as was recently shown
by Bryant~\cite{Br}, the set of missing, or {\em exotic}, holonomies is
non-empty. As usual in the representation theory, in order to get a deeper
understanding of all irreducible real holonomies one should first
try to address a complex version of the problem. The main result announced
in this paper asserts that any torsion-free holomorphic affine connection
with irreducibly acting holonomy group can be generated by twistor methods.

\paragraph{1.\ Irreducible $G$-structures.}
Let $M$ be an $m$-dimensional complex manifold and ${\cal L}^*M$ the
holomorphic coframe bundle $\pi: {\cal L}^*M \rar M$ whose fibres
${\cal L}^*_tM =\pi^{-1}(t)$ consist of all $\Bbb C$-linear isomorphisms
$e: \Bbb C^{m} \rar \Omega^1_t M$. The space ${\cal L}^*M$ is a principle
right $GL(m,\Bbb C)$-bundle with the right action given by $R_g(e)= e\circ g$.
If $G$ is a closed subgroup of $GL(m,\Bbb C)$, then a (holomorphic)
$G$-structure on $M$ is a principle subbundle $\cal G$ of ${\cal L}^*M$ with
the group $G$. It is clear that there is a one-to-one correspondence between
the set of $G$-structures on $M$ and holomorphic sections $\sigma$ of the
quotient bundle $\tilde{\pi}: {\cal L}^*M/G \rar M$ whose typical fibre is
isomorphic to $GL(m,\Bbb C)/ G$. A $G$-structure on $M$ is called {\em locally
flat}\, if ${\cal L}^*M/G$ can be trivialised over a sufficiently small
neighbourhood, $U$, of each point $t\in M$ in such a way that the associated
section $\sigma$ of ${\cal L}^*M/G$ is represented over $U$ by a constant
$GL(m,\Bbb C)/G$-valued function. A G-structure is called {\em k-flat}\, if,
for each $t\in M$, the $k$-jet of the associated section $\sigma$ of
${\cal L}^*M/G$ at $t$ is isomorphic to the $k$-jet of some locally flat
section of ${\cal L}^*M/G$. It is easy to show that a $G$-structure admits
a torsion-free affine connection if and only if it is 1-flat (cf.\
\cite{Bryant2}). A $G$-structure on $M$ is called {\em irreducible}\, if the
action of $G$ on $\Bbb C^{m}$ leaves no non zero invariant subspaces.

When studying a (torsion-free) affine connection $\nabla$ on a connected
simply connected complex manifold $M$ with the irreducibly acting holonomy
group $G$, one usually works with the associated irreducible (1-flat)
$G$-structure $\cal G_{\nabla}\subset{\cal L}^*M$ \cite{Bryant1,S}. Define
two points $u$
and $v$ of ${\cal L}^*M$ to be equivalent, $u\sim v$, if there is a
holomorphic path $\gamma$ in $M$ from $\pi(u)$ to $\pi(v)$ such that
$u=P_{\gamma}(v)$, where $P_{\gamma}:\Omega^1_{\pi(v)}M\rar\Omega^1_{\pi(u)}M$
is the parallel transport along $\gamma$. Then $\cal G_{\nabla}$ can be
defined, up to an isomorphism, as $\left\{u\in {\cal L}^*M \mid u\sim v
\right\}$ for some coframe $v$. The $G$-structure ${\cal G}_{\nabla}$ is
the smallest subbundle of ${\cal L}^*M$ which is invariant under
$\nabla$-parallel translations.

\paragraph{2.\ Complex contact structures.}
Let $Y$ be a complex $(2n + 1)$-dimensional manifold. A complex contact
structure on $Y$ is a rank $2n$ holomorphic subbundle $D\subset TY$ of the
holomorphic tangent bundle to $Y$ such that the Frobenius form
\begin{eqnarray*}
\Phi: D \times D & \longrightarrow & TY/D\\
(v,w) & \longrightarrow & [v,w]\bmod D
\end{eqnarray*}
is non-degenerate.
A complex $n$-dimensional submanifold $X$ of the complex contact manifold $Y$
is called a {\em Legendre submanifold}\, if $TX\subset D$. The normal bundle
of a Legendre submanifold $X\hookrightarrow Y$ is isomorphic to $J^1 L_X$
\cite{L2}, where $L_X = \left.L\right|_X$ and $L$ is the contact line bundle
on $Y$ defined by the exact sequence
$$
0\lon D \lon TY  {\lon} L \lon 0.
$$
Given a Legendre submanifold
$X\hook Y$, there is a naturally associated "flat" model, $X\hook J^1 L_X$,
consisting of the total space of the vector bundle
$J^1 L_X$ together with its canonical contact structure
and the Legendre submanifold $X$ realised as a zero section of
$J^1 L_X \rar X$. The Legendre submanifold $X\hook Y$ is called
$k$-{\em flat}\, if the $k$\,th-order Legendre jet \cite{L2} of $X$ in $Y$
is isomorphic to the $k$\,th-order Legendre jet of $X$
in $J^1 L_X$. The obstruction for a complex Legendre submanifold to be
1-flat is a  cohomology class in
$H^1\left(X, L_X\otimes S^2(J^1 L_X)^*\right)$.

\paragraph{3.\ Twistor theory and $G$-structures.}
Recall that a generalised flag variety $X$ is a compact
simply connected homogeneous K\"{a}hler manifold \cite{BE}. Any such a
manifold is of the form $X=G/P$, where $G$ is a complex semisimple Lie
group and $P\subset G$ a fixed parabolic subgroup.

\begin{theorem}\label{Main}
Let $X$ be a generalised flag variety embedded as a Legendre submanifold
into a complex contact manifold $Y$ with contact line bundle $L$ such that
$L_X$ is very ample on $X$. Then
\begin{description}
\item[(i)] There exists a maximal family $\{X_t \hook Y\mid t\in M\}$ of
compact complex Legendre submanifolds obtained by holomorphic deformations
of $X$ inside $Y$. Each submanifold $X_t$ is isomorphic to $X$. The moduli
space $M$, called a Legendre moduli space, is an $m$-dimensional complex
manifold, where $m=h^0(X,L_X)$.

\item[(ii)] The Legendre submanifold $X\hook Y$ is stable under holomorphic
deformations of the contact structure on (the tubular neighbourhood of
$X$ in) $Y$.

\item[(iii)] For each $t\in M$, there is a canonical isomorphism
$s: T_t M\rightarrow H^0(X_t, L_{X_t})$ representing a tangent vector at $t$
as a global holomorphic section of the line bundle
$L_{X_t} = \left.L\right|_{X_t}$.

\item[(iv)] The Legendre moduli space $M$ comes equipped with an induced
irreducible $G$-structure, ${\cal G}_{ind} \rar M$,
with $G$ isomorphic to the connected component of the identity of
the group of all global biholomorphisms $\phi: L_X
\rightarrow L_X$ which commute with the projection $\pi: L_X\rightarrow X$.

\item[(v)]\ The induced $G$-structure $\cal G_{ind}$ is 1-flat
(i.e.\ torsion-free) if and only if the complete family
$\{X_t\hook Y\mid t\in M\}$ consists of 1-flat Legendre submanifolds.

\item[(vi)] If ${\cal G}_{ind}$ is 1-flat, then the bundle of all torsion-free
connections in ${\cal G}_{ind}$ has as the typical fiber an affine space
modelled on $H^0\left(X, L_{X}\otimes S^2(J^1 L_{X})^*\right)$.

\item[(vii)] If ${\cal G}_{ind}$ is 1-flat, then the obstruction for
${\cal G}_{ind}$ to be 2-flat is given by a tensor field on $M$
whose value at each $t\in M$ is represented by a cohomology class
$\rho_t\in H^1\left(X_t, L_{X_t}\otimes S^3(J^1 L_{X_t})^*\right)$.

\item[(viii)] Let $H\subset GL(k, \Bbb C)$ be one of the following
subgroups: (a) $SO(2n+1,\Bbb C)$ when $k=2n+2\geq 8$; (b) $Sp(2n+2,\Bbb C)$
when $k=2n+2\geq 4$; (c) $G_2$ when $k=7$. Suppose that $G\subset GL(m,\Bbb C)$
is a connected semisimple Lie subgroup whose decomposition into a
locally direct
product of simple groups contains $H$. If $\cal G$ is any irreducible
1-flat $G\times \Bbb C^*$-structure on  an $m$-dimensional manifold $M$,
then there exists a complex contact manifold $(Y,L)$ and a generalised flag
variety $X$ embedded into $Y$ as a Legendre submanifold with $L_X$ being
very ample, such that, at least locally, $M$ is canonically isomorphic to the
associated Legendre moduli space and $\cal G \subset \cal G_{ind}$. In
particular, when $G=H$ one has in the case
(a) $X= SO(2n+2,\Bbb C)/U(n+1)$ and $\cal G_{ind}$ is a
$CO(2n+2,\Bbb C)$-structure; in the case (b)
$X=\Bbb C \Bbb P^{2n+1}$ and $\cal G_{ind}$ is a $GL(2n+2,\Bbb C)$-structure;
and in the case (c) $X=Q_5$ and $\cal G_{ind}$ is a
$CO(7,\Bbb C)$-structure.

\item[(ix)] Let $G\subset GL(m,\Bbb C)$ be an arbitrary connected
semisimple Lie
subgroup whose decomposition into a locally direct product of
simple groups does not contain any of the groups $H$ considered in (viii).
If $\cal G$ is any irreducible 1-flat $G\times \Bbb C^*$-structure on an
$m$-dimensional manifold $M$, then there exists a complex contact manifold
$(Y,L)$ and a Legendre submanifold $X\hook Y$ with $X=G/P$ for some parabolic
subgroup $P\subset G$ and with $L_X$ being very ample, such that, at least
locally, $M$ is canonically isomorphic to the associated Legendre moduli
space and $\cal G = \cal G_{ind}$.
\end{description}
\end{theorem}

The Lie algebra of the group $G$ of all global biholomorphisms $L_X
\rar L_X$ which commute with the projection $\pi: L_X\rightarrow X$ is exactly
the vector space $H^0\left(X, L_X\otimes (J^1 L_X)^*\right)$ with its
natural Lie algebra structure \cite{Me1}. If $X=G/P$, then $\cal G_{ind}$
on the associated Legendre moduli space is often isomorphic to a
$G\times {\Bbb C}^*$-structure, but there are exceptions \cite{Ahiezer}
which are considered in Theorem~\ref{Main}(viii). In these exceptional cases
the original $G\times \Bbb C^*$-structure may not be equal to the induced one,
and one might try to identify some additional structures on the associated
twistor spaces $(Y,L)$ which ensure that $\cal G_{ind}$ admits a necessary
reduction. However, the "exceptional" $H$-structures with $H$ as in
(a), (b) and (c) of Theorem~\ref{Main}(viii) are fairly well
understood by now \cite{Berger,Bryant1,Bryant2,S}. If there exists an exotic
torsion-free $G$- or $G\times \Bbb C^*$-structure with simple $G$ other than
Bryant's $H_3$ \cite{Br}, it must be covered, up to a $\Bbb C^*$ action, by
the "generic" clause (ix) in Theorem~\ref{Main}.

Two particular examples of this general construction have been
considered earlier \cite{L1,Br}.
The first example is a pair $X\hook Y$ consisting of an $n$-quadric $Q_n$
embedded into a $(2n+1)$-dimensional contact manifold $(Y,L)$ with
$\left. L\right|_X\simeq i^*\f_{{\Bbb C \Bbb P}^{n+1}}(1)$,
$i: Q_n\hook {\Bbb C \Bbb P}^{n+1}$ being a standard projective realisation of
$Q_n$. It is easy to check that in this case \mbox{$H^0(X,L_X\otimes
(J^1 L_X)^{*})$} is precisely the conformal algebra implying that the
associated $(n+2)$-dimensional Legendre moduli space $M$ comes equipped
canonically with a conformal structure. This is in accord with LeBrun's paper
\cite{L1}, where it has been shown how a conformal Weyl connection can be
encoded into a complex contact structure on the space of complex null
geodesics.
Since $H^1\left(X, L_X\otimes S^2(J^1 L_X)^*\right)=0$, the induced
conformal structure must be torsion-free in agreement with the classical
result of differential geometry. Easy calculations show that the vector space
$H^1\left(X, L_X\otimes S^3(J^1 L_X)^*\right)$ is exactly the subspace
of $TM\otimes\OM \otimes \Omega^2M$ consisting of tensors with
Weyl curvature symmetries. Thus Theorem~\ref{Main}(vii) implies the well-known
Schouten conformal flatness criterion. Since $H^0\left(X, L_X\otimes S^2
(J^1 L_X)^*\right)$ is isomorphic to the typical fibre of $\OM$, the set
of all torsion-free affine connections preserving the induced conformal
structure is an affine space modelled on $H^0(M,\OM)$, again in agreement
with the classical result.

The second example, which also was among motivations behind the present work,
is Bryant's \cite{Br} relative deformation problem $X\hookrightarrow Y$ with
$X$ being a rational Legendre curve ${\Bbb C\Bbb P}^1$ in a complex contact
3-fold $(Y,L)$ with $L_X=\f (3)$. Calculating
$H^0(X,L_X\otimes (J^1 L_X)^{*})$, one easily concludes that the induced
$G$-structure, $\cal G_{ind}$, on the associated
4-dimensional Legendre moduli space is exactly an exotic $G_3$-structure
which has been studied by Bryant in his search for irreducibly acting
holonomy groups of torsion-free affine connections which are missing
in the Berger list \cite{Berger}. Since
$H^1\left(X, L_X\otimes S^2(J^1 L_X)^*\right)=0$, Theorem~\ref{Main}(v) says
 the induced $G_3$-structure $\cal G_{ind}$ is torsion-free in accordance
with \cite{Br}. Since $H^0\left(X, L_X\otimes S^2(J^1 L_X)^*\right)=0$,
$\cal G_{ind}$ admits a unique torsion-free affine connection $\nabla$.
The cohomology class $\rho_t\in H^1\left(X, L_X\otimes S^3(J^1 L_X)^*\right)$
from Theorem~\ref{Main}(v) is exactly the curvature tensor of $\nabla$.

\paragraph{4.\ Outline of the proof of Theorem~\ref{Main}.}
Items (i)-(iii) of Theorem~\ref{Main} follow from a more general theorem
\cite{Me1} which says that if $X$ is a compact complex Legendre submanifold
of a complex contact manifold $(Y,L)$ such that $H^{1}(X,L_{X}) = 0$, then
there exists a complete, maximal and stable analytic family
$\{X_{t}\hook Y\mid t\in M\}$ of compact Legendre submanifolds containing
$X$ (completeness, e.g., means that the natural map $T_t M\rar
H^0(X_t, L_{X_t})$ is isomorphic for each $t\in M$). Indeed, if $X$ is a
generalised flag variety and $L_X$ a very ample line bundle on $X$, then
$h^0(X,L_X)> 0$, and hence $H^1(X,L_X)=0$ by Bott-Borel-Weil theorem and
the fact that any holomorphic line bundle on $X$ is homogeneous. Since $X$
is rigid, each submanifold $X_t$ of the family is isomorphic to $X$.

In view of the canonical isomorphism $T_t M\rar H^0(X_t, L_{X_t})$
the item (iv) is not a surprise. More precisely, defining
$F=\{(y,t)\in Y\times M \mid y\in X_t\}$ and using the fact that
$L_{X_t}$ is very ample on $X_t$, one can easily realise $F$ as a
subbundle of the projectivised conormal bundle $\Bbb P_M(\OM)$. Fibrewise,
this construction is the well-known projective realisation of a
generalised flag variety $X$ in
${\Bbb C\Bbb P}^{m-1}\simeq \Bbb P\left(H^0(X,L_X)^*\right)$ \cite{BE}.
The subgroup $G$ of $GL(m,\Bbb C)$ which leaves $X\subset \Bbb C\Bbb P^{m-1}$
invariant is exactly the one described in item (iv) of Theorem~\ref{Main}.

The next question we address is  how to distinguish in terms of the holomorphic
embedding data $X\hook Y$ the subclass of 1-flat induced $G$-structures.
This leads us to explore the towers of infinitesimal neighbourhoods of two
embeddings of analytic spaces, $X_t\hook Y$ and $t\hook M$. At the first
floors of these towers we have, by item (iii) of Theorem~\ref{Main}, an
isomorphism $T_tM = H^0(X_t, L_{X_t})$ which is in the basis of the above
conclusion about the induced $G$-structure on $M$. The second floors of these
two towers are related to each other as follows.
If $J_t\subset \f_M$ is the ideal of holomorphic functions which vanish
at $t\in M$, then the tangent space $T_tM$   is isomorphic to
$\left(J_t/J_t^2\right)^*$. Define a second order tangent bundle,
$T_t^{[2]}M$, at the point $t$ as $\left(J_t/J_t^3\right)^*$. This definition
implies that $T_t^{[2]}M$ fits into an exact sequence of complex vector spaces
\begin{equation}
0\lon T_tM \lon T_t^{[2]}M \lon S^2(T_tM) \lon 0 \label{mu}
\end{equation}
For each $t\in M$ there
exists a holomorphic line bundle, $\Delta_{X_t}^{[2]}$, on the
associated Legendre submanifold $X_t\hook Y$ such that there are an exact
sequence of locally free sheaves
\begin{equation}
0\lon L_{X_t} \stackrel{\alpha}{\lon} \Delta_{X_t}^{[2]} \lon S^2(N_t) \lon 0
\label{mumu}
\end{equation}
and a commutative diagram of vector spaces
\begin{equation}
\begin{array}{rccccccccl}
0 & \lon & T_t M & \lon & T^{[2]}_t M & \lon & S^2(T_t M) & \lon & 0 \\
&& \Big\downarrow && \Big\downarrow && \Big\downarrow && \\
0 & \lon & H^0(X_t,L_{X_t}) & \lon & H^0\left(X_t,\Delta_{X_t}^{[2]}\right) &
\lon & H^0\left(X_t, S^2(N_{X_t})\right) & \lon & 0
\end{array} \label{mumumu}
\end{equation}
which extends the canonical isomorphism $T_tM \rightarrow H^0(X_t, L_{X_t})$ to
{\em second}\, order infinitesimal neighbourhoods of $t\hookrightarrow M$
and $X_t\hookrightarrow Y$. All we need to know in this paper about
$\Delta_{X_t}^{[2]}$ is that this bundle exists and has the stated
properties. For details of its definition we refer the interested reader
to \cite{Me1,Me2}.

One can show that the Legendre submanifold $X_t\hook Y$ is 1-flat if and
only if the obstruction, $\delta_{X_t}^{[2]}\in
H^1\left(X_t, L_X\otimes S^2(J^1 L_{X_t})^*\right)$, for the global
splitting of (\ref{mumu}) is zero. If $X_t$ is 1-flat, then any splitting
of (\ref{mumu}), i.e.\ a morphism $\beta: \Delta_{X_t}^{[2]}
\rightarrow L_{X_t}$ such that $\beta\circ\alpha = id$, induces via the above
commutative diagram an associated splitting of the exact sequence
(\ref{mu}) which is  equivalent to a torsion-free affine connection at
$t\in M$. This implies that if the family $\{X_{t}\hook Y\mid t\in M\}$
consists of 1-flat Legendre submanifolds, then the induced $G$-structure
admits a torsion-free affine connection, i.e.\ is 1-flat. In reverse order,
given a Legendre moduli space $M$ with the induced $G$-structure
$\cal G_{ind}$ as in Theorem~\ref{Main}(iv), one can use the commutative
diagram~(\ref{mumumu}) to show that that any torsion-free connection
in $\cal G_{ind}$ induces canonically a splitting of the exact
sequence~(\ref{mumu}) for each $t\in M$. The set of all splittings
of~(\ref{mumu}) is an affine space modelled on
$H^0\left(X_t, L_{X_t}\otimes S^2(N_t^*)\right)\simeq
H^0\left(X, L_{X}\otimes S^2(N^*)\right)$. These facts prove items (v)
and (vi) of Theorem~\ref{Main}. The item (vii) can be proved by
straightforward calculations in Darboux local coordinates on $Y$. For details
of these calculations we refer to \cite{Me1}.

Consider now a complex $m$-dimensional manifold $M$ and an
irreducible $G$- or $G\times \Bbb C^*$-structure $\cal G\subset \cal L^*M$,
where $G\subset GL(m,\Bbb C)$ is a semisimple Lie subgroup (any irreducible
$H$-structure with reductive $H$ must be one of these). Since $\cal G$ is
irreducible, there is a naturally associated subbundle $\tilde{F}\subset
\OM$ whose typical fiber is the cone in $\Bbb C^m$ defined as
the $G$-orbit of the line spanned by a highest weight vector.
Denote $\tilde{\cal F} = \tilde{F}\setminus 0_{\tilde{F}}$, where
$0_{\tilde{F}}$ is the "zero" section of $\tilde{p}: \tilde{F} \rar M$ whose
value at each  $t\in M$ is the vertex of the cone $\tilde{p}^{-1}(t)$.
The quotient bundle $\bnu: {\cal F} \equiv \tilde{\cal F}/{\Bbb C^*} \lon M$
is then a subbundle of the projectivized cotangent bundle $\Bbb P_M(\OM)$
whose fibres $X_t$ are isomorphic to the generalised flag variety
$G/P$, where $P$ is the parabolic subgroup of $G$ which preserves the
highest weight vector in $\Bbb C^m$ up to a scale. Denote $\dim G/P = n$.
The total space of the cotangent bundle $\OM$ has a canonical holomorphic
symplectic 2-form $\omega$. Then the sheaf of holomorphic functions on
$\OM$ is a sheaf of Lie algebras relative to the Poisson bracket
$\{f,g\}= \omega^{-1}(df,dg)$.
The pullback, $i^*{\omega}$, of the symplectic form $\omega$ from
$\OM\setminus 0_{\OM}$ to its submanifold $i:
\tilde{\cal F}\lon \OM \setminus 0_{\OM}$ defines a distribution
${\cal D}\subset T\tilde{\cal F}$ as the kernel of the natural "lowering of
indices" map $T\tilde{\cal F} \stackrel{\lrcorner\,\omega}{\lon}
\Omega^1\tilde{\cal F}$, i.e.\
${\cal D}_e = \left\{V\in T_e\tilde{\cal F}: V\lrcorner\, i^*\omega =
0\right\}$
at each point $e\in \tilde{\cal F}$. Using the fact that $d(i^*\omega) =
i^*d\omega = 0$, one can show that this distribution is integrable and thus
defines  a foliation of $\tilde{\cal F}$ by holomorphic leaves.
We shall assume that the space of leaves, $\tilde{Y}$, is a complex manifold.
This assumption imposes no restrictions on the local structure of $M$.
The fact that the Lie derivative, ${\cal L}_V i^*\omega =
V\lrcorner\, i*d\omega + d(V\lrcorner\, i^*\omega) = 0$, vanishes for any
vector field
$V$ tangent to the leaves implies that $i^*\omega$ is the pullback
relative to the canonical projection $\tilde{\bmu}: \tilde{\cal F}
\rar \tilde{Y}$ of a closed 2-form $\tilde{\omega}$ on $\tilde{Y}$. It is
easy to check that $\tilde{\omega}$ is non-degenerate. The quotient
$(\tilde{Y}, \tilde
{\omega})$ is what is usually called a symplectic reduction of
$(\OM\setminus 0_{\OM}, \omega)$ along the submanifold $\tilde{\cal F}$
 (cf.\ \cite{S}).

Let $e$ be any point of $\tilde{\cal F}\subset \OM\setminus 0_{\OM} $.
Restricting a "lowering of indices" map $T_e(\OM) \stackrel{\lrcorner\,\omega}
{\lon}\Omega_e(\OM)$ to the subspace ${\cal D}_e$, one obtains an injective map
$$
0\lon {\cal D}_e\stackrel{\lrcorner\,\omega}{\lon} {\cal N}_e^*,
$$
where ${\cal N}_e^*$ is the fibre of the conormal bundle of
$\tilde{\cal F}\hook\OM\setminus 0_{\OM}$. Therefore, the rank of the
distribution ${\cal D}$ is equal at most to $\rank\, {\cal N}^* =m-n-1$.
It is easy to show that $\rank\, \cal D$ is maximal possible if and only if
the ideal
sheaf, $I_{\tilde{\cal F}}$, of $\tilde{\cal F}$ in $\OM\setminus 0_{\OM}$
is a sheaf of Lie subalgebras, i.e.\ $\left\{I_{\tilde{\cal F}},
I_{\tilde{\cal F}}\right\}\subset I_{\tilde{\cal F}}$. An irreducible
$G$-structure is called {\em Poisson}\, if $I_{\tilde{\cal F}}$ is a Lie
subalgebra (equivalently, if $\rank\, \cal D = m-n-1$). It is easy to check
that
a locally flat $G$-structure is Poisson. Since $\{\ ,\ \}$ is a first order
differential operation, this immediately implies that  any irreducible
1-flat $G$-structure is also Poisson.

Next we show that if $G$ is a reductive Lie group, then every complex
$m$-manifold $M$ with a Poisson $G$-structure is canonically
isomorphic, at least locally, to a Legendre moduli space. Indeed, there is
an integrable distribution $\cal D$ of rank $m-n-1$ on the bundle
$\tilde{\cal F}\rar M$. Assuming that $M$ is sufficiently "small", we obtain
as the quotient a symplectic manifold, $(\tilde{Y},\tilde{w})$, with
$\dim \tilde{Y}= (m+n+1)-(m-n-1) = 2n+2$. There is a double fibration
$$
\tilde{Y} \stackrel{\tilde{\bmu}}{\longleftarrow} \tilde{\cal F}
\stackrel{\tilde{\bnu}}{\lon} M
$$
with fibres of $\tilde{\bmu}$ being leaves of the integrable distribution
$\cal D$ and fibres of $\tilde{\bnu}$ being $(n+1)$-dimensional cones
$X_t\subset \Omega^1_t M\setminus 0$ generated by $G$-orbits of highest weight
vectors. It is clear that, for each
$t\in M$, the submanifold $\tilde{\bmu}(X_t)\subset \tilde{Y}$ is a Lagrange
submanifold relative to the induced symplectic form $\tilde{\omega}$ on
$\tilde{Y}$.

There is a natural action of $\Bbb C^*$ on $\tilde{\cal F}$
which leaves ${\cal D}$ invariant and thus induces an action of $\Bbb C^*$ on
$\tilde{Y}$. The quotient $Y = \tilde{Y}/\Bbb C^*$ is a $(2n+1)$-dimensional
complex manifold which has a double fibration structure
$$
Y \stackrel{\bmu}{\longleftarrow} {\cal F}=\tilde{\cal F}/\Bbb C^*
\stackrel{\bnu}{\lon} M
$$
and thus contains a family of compact $n$-dimensional embedded
manifolds\linebreak
$\left\{X_t =\bmu\circ\bnu^{-1}(t)\hook Y\mid t\in M\right\}$ with
$X_t = \tilde{X}_t/\Bbb C^*$. Next, inverting a well-known procedure of
symplectivisation of a contact manifold, it is not hard to show
that $Y$ has a complex contact structure such that
the family $\left\{X_t\hook Y\mid t\in M\right\}$ is a family of
compact Legendre submanifolds. The contact line bundle $L$ on $Y$
is just the quotient $L=\tilde{\cal F}\times  \Bbb C/\Bbb C^*$ relative to the
natural multiplication map $\tilde{\cal F}\times\Bbb C \lon
\tilde{\cal F}\times  \Bbb C$, $(p,c)\rar (\lambda p, \lambda c)$, where
$\lambda\in \Bbb C^*$. Then $\left.L\right|_{X_t}$ is isomorphic to
the restriction of the hyperplane section bundle on $\Bbb P(\Omega^1_t M)$
to its submanifold $\bnu^{-1}(t)\simeq X_t$ and hence is very ample on
$X_t$. Therefore, $h^0\left(X_t,\left.L\right|_{X_t}\right) = m$ which
implies that the Legendre family $\left\{X_t\hook Y\mid t\in M\right\}$ is
complete.

Therefore we proved that if $G\subset GL(m,\Bbb C)$ is a semisimple Lie
subgroup and $\cal G$ any irreducible 1-flat $G\times \Bbb C^*$-structure on
an $m$-dimensional manifold $M$, then
there exists a complex contact manifold $(Y,L)$ and a Legendre submanifold
$X\hook Y$ with $X=G/P$ for some parabolic subgroup $P\subset G$ and with
$L_X$ being
very ample, such that, at least locally, $M$ is canonically isomorphic to the
associated Legendre moduli space. To complete the proof of final items (viii)
and (ix) of Theorem~\ref{Main} one needs only the fact \cite{Ahiezer} that
if $X$ is an irreducible
generalised flag variety of a connected complex semisimple Lie
group $G$, then the connected component, $\mbox{Aut}^0{X}$, of the group of
global holomorphic automorphisms is simple and, as a rule, coincides with
$G$. The only exceptions are listed below:
(a) $X= SO(2n+2,\Bbb C)/U(n+1)$, $n\geq 3$, $G=SO(2n+1,\Bbb C)$,
$\mbox{Aut}^0{X}= SO(2n+2,\Bbb C)$; (b) $X=\Bbb C \Bbb P^{2n+1}$, $n\geq 1$,
$G=Sp(2n+2, \Bbb C)$, $\mbox{Aut}^0{X} = SL(2n+2,\Bbb C)$; (c) $X=Q_5$, the
5-dimensional compact quadric, $G=G_2$, $\mbox{Aut}^0{X}=SO(7,\Bbb C)$.
This fact completes the outline of the proof of Theorem~\ref{Main}. $\Box$

\paragraph{5. On holonomy groups.} Let $\{X_{t}\hook Y\mid t\in M\}$ be a
complete family of compact Legendre submanifolds. A torsion-free connection
on $M$ which arises at each $t\in M$ from a splitting of the
extension~(\ref{mumu}) is called an {\em induced}\, connection. In the
previous subsection we proved also the following
\begin{theorem}
Let $\nabla$ be a holomorphic torsion-free affine connection on a complex
manifold $M$ with irreducibly acting reductive holonomy group $G$. Then
there exists a complex contact manifold $(Y,L)$ and a Legendre submanifold
$X\hook Y$ with $X=G_s/P$ for some parabolic subgroup $P$ of the semisimple
factor $G_s$ of $G$ and with $L_X$ being very ample, such that, at least
locally, $M$ is canonically isomorphic to the associated Legendre moduli space
and $\nabla$ is an induced torsion-free affine connection in $\cal G_{ind}$.
\end{theorem}

Theorem 2 and much of Theorem~1 are devoted to torsion-free affine
connections and 1-flat  $G$-structures. In fact, as follows from the outline
of the proof of Theorem~1, the class of irreducible $G$-structures that can
be interpreted as induced on Legendre moduli spaces is
much larger than the class of 1-flat $G$-structures (but much smaller
than the class of all possible irreducible $G$-structures). For motives
explained in sect.4, $G$-structures in this class are called Poisson.
The theorem of Hano and Ozeki \cite{HO} is no longer true in the category
of affine connections in Poisson $G$-structures which implies that,
in addition to the open problem of classifying all irreducibly acting
holonomies of torsion-free affine connections, we get another
seemingly non-trivial problem of classifying all irreducibly acting
holonomies $G$ of affine connections with {\em non-zero torsion}\, which
are tangent to Poisson $G$-structures.

We conclude this paper with the remark that given a sufficiently
"small" complex $m$-dimensional manifold $M$ and an irreducible
$G$-structure $\cal G\subset \cal L^*M$ with reductive $G$, one may proceed
as in sect.4 to show that there is an associated complex contact
$(2q+1)$-dimensional manifold $(Y,D)$ and $m$-dimensional family of
generalised flag varieties $\{X_t\hook Y\mid t\in M\}$ which are tangent
to the contact distribution $D$, i.e.\ $TX_t\subset D$ for each $t\in M$.
If $\cal G$ is 1-flat or, more generally, Poisson, then each contact
submanifold $X_t$ is of maximal possible dimension $q$, i.e.\ it is a
Legendre submanifold. In general, $\dim X_t\leq q$, and we
can stratify all possible $G$-structures on $M$ into classes parametrised
by an integer $l= q - \dim X_t$. In this paper it is shown that in the case
$l=0$ (much of) the original $G$-structure together with its
basic geometric invariants can be reconstructed from the
complex contact structure on $Y$ by twistor methods. It seems unlikely
that something similar can be done in the case $l\geq 1$.

\vspace{3 mm}

\noindent{\em Acknowledgements}.
It is a pleasure to thank D.V.\ Alekseevski, T.\ Bailey, M.\  Eastwood,
S.\ Huggett, L.\ Mason, H.\ Pedersen, Y.-S.\ Poon and P.\ Tod for valuable
comments and remarks.
I am  grateful to the University of Odense, St.\ John's College of Oxford
University and the Erwin Schr\"{o}dinger International Institute for
hospitality during the work on different stages of this project.


\begin{thebibliography}{99}

\bibitem[A]{Ahiezer} {\sc D.N.\ Ahiezer}, {\em Homogeneous complex
manifolds}, in {\em Several Complex Variables IV}, Springer 1990.

\bibitem[Al]{Al} {\sc D.V.\ Alekseevski}, 'Riemannian spaces with unusual
holonomy groups', {\em Funct.\ Anal.\ Appl.}\, {f 2} (1968), 97-105.

\bibitem[B-E]{BE} {\sc R.J.\ Baston and M.G.\ Eastwood}, {\em The Penrose
transform, its interaction with representation theory} (Oxford University
Press, 1989).

\bibitem[B]{Berger} {\sc M.\ Berger}, 'Sur les groupes d'holonomie des
vari\'{e}t\'{e}s \'{a} connexion affine et des vari\'{e}t\'{e}s
Riemanniennes',  {\em Bull.\ Soc.\ Math.\ France}\, { 83} (1955), 279-330.

\bibitem[Br1]{Bryant1} {\sc R.\ Bryant}, 'A survey of Riemannian metrics
with special holonomy groups', In: {\em Proceedings of the ICM at Berkley},
1986. Amer.\ Math.\ Soc., 1987, 505-514.

\bibitem[Br2]{Bryant2} {\sc R.\ Bryant}, 'Metrics with exceptional
holonomy', {\em Ann.\ of Math.}\, (2) { 126} (1987), 525-576.

\bibitem[Br3]{Br} {\sc R.\ Bryant}, 'Two exotic holonomies in dimension
four, path geometries, and twistor theory', {\em Proc.\ Symposia in Pure
Mathematics}\, { 83} (1991), 33-88.

\bibitem[H-O]{HO} {\sc J.\ Hano and H.\ Ozeki}, 'On the holonomy groups
of linear connections' {\em Nagoya Math.\ J.}\, { 10} (1956), 97-100.

\bibitem[K]{Kodaira} {\sc K.\ Kodaira}, 'A theorem of completeness of
characteristic systems for analytic families of compact submanifolds of
complex manifolds' {\em Ann.\ Math.}\, { 75} (1962), 146-162.

\bibitem[L1]{L1} {\sc C.\ R.\ LeBrun}, 'Spaces of complex null geodesics
in complex-Riemannian geometry' {\em Trans.\ Amer.\ Math.\ Soc.}\, { 284}
(1983), 209-321.

\bibitem[L2]{L2} {\sc C.\ R.\ LeBrun}, 'Thickenings and conformal gravity'
{\em Commun.\ Math.\ Phys.}\, { 139} (1991), 1-43.

\bibitem[Me1]{Me1} {\sc S.\ A.\ Merkulov}, 'Existence and geometry
of Legendre moduli spaces', preprint.

\bibitem[Me2]{Me2} {\sc S.\ A.\ Merkulov}, 'Moduli of compact complex
Legendre submanifolds of complex contact manifolds', To appear
in {\em Math. Research Lett.}.

\bibitem[P]{P} {\sc R.\ Penrose}, 'Non-linear gravitons and curved twistor
theory' {\em Gen.\ Rel.\ Grav.}  { 7} (1976), 31-52.


\bibitem[S]{S} {\sc S.\ M.\ Salamon}, {\em Riemannian Geometry and
Holonomy Groups} (Longman, 1989).

\end{thebibliography}
\end{document}